\documentclass[12pt,a4paper, twoside, reqno]{amsart} 

\usepackage[margin=3cm]{geometry}

\usepackage{ifpdf}
\usepackage{graphicx}
\usepackage[T1]{fontenc} 
\usepackage[latin1]{inputenc} 
\usepackage{amsmath} 
\usepackage{amsthm} 
\usepackage{amsfonts} 
\usepackage{amssymb}
\usepackage{amsbsy}
\usepackage{bbm} 
\usepackage[english]{babel} 
\usepackage{varioref} 
\usepackage{enumerate} 
 \usepackage{setspace}
\usepackage{fancyhdr}

\theoremstyle{definition}

\newtheorem*{defn*}{Definition}

\theoremstyle{remark}

\theoremstyle{plain}
\newtheorem{thm}{Theorem}[section]
\newtheorem*{thm*}{Theorem}
\newtheorem{lem}[thm]{Lemma}

\newtheorem{cor}[thm]{Corollary}

\newtheorem*{thmharrap}{Theorem A (2011)}

\newcommand{\norm}[1]{\ensuremath{\left\Vert #1 \right\Vert}}

\newcommand{\infabs}[1]{\ensuremath{\left\vert #1 \right\vert}}

\newcommand{\Bad}{\text{\textup{\bf Bad}}}

 \DeclareMathOperator{\mat}{Mat}

\numberwithin{equation}{section}

\title[A note on weighted BA affine forms]{A note on weighted badly approximable linear forms}

\author{Stephen Harrap$^\dag$}
\thanks{$^\dag$ Research supported by EPSRC grant number $EP/L005204/1$l}

\address{S. Harrap, Durham University, Department of Mathematical Sciences, Science Laboratories, South Rd, Durham,
	DH1 3LE, United Kingdom}
\email{s.g.harrap@durham.ac.uk}

\author{Nikolay Moshchevitin$^*$}
 \thanks{ $^*$ Research supported by RNF grant 14-11-00433}
 
\address{N. Moshchevitin, Leninskie Gory 1, GZ MGU, Department of Mathematics and Mechanics, Moscow State University,
  119991 Moscow, Russia}
\email{moshchevitin@rambler.ru}

\date{}

\begin{document}
	
	\subjclass[2010]{Primary: 11J83; Secondary 11J13, 11K60}
	\keywords{Metric Diophantine approximation; twisted inhomogeneous; circle rotations; Schmidt games, winning sets}

\begin{abstract}
We prove a result in the area of twisted Diophantine approximation related to the theory of Schmidt games. In particular, under certain restrictions we give a affirmative answer to the analogue in this setting of a famous conjecture of Schmidt from Diophantine approximation. 
\end{abstract}

\maketitle
\section{Introduction}
\label{sec:introduction}

In 2007, Kim \cite{Kim} proved that for any irrational number $x$ the set of real $\alpha \in [0,1)$ for which the property
\begin{equation}
	\label{eqn:Kim}
	\inf_{q \in \mathbb{N}} q \norm{qx - \alpha } \: = \: 0
\end{equation}
holds has maximal Lebesgue measure $1$. Here and throughout,  $\norm{ \, . \, }$ denotes the distance to the nearest integer. 
This statement has a tangible interpretation in terms of a rotation of the unit circle. Identifying the unit circle with the unit interval $[0,1)$, the value $qx$ (modulo $1$) may be thought of as the position of the origin after $q$ rotations by the angle $x$. A celebrated result of Weyl~\cite{MR1511862} implies that every such irrational rotation in $[0,1)$ visits any fixed set of positive measure infinitely often. To be precise, Kim's result corresponds to the situation when this set of positive measure (realized here by an interval centred at the point $\alpha$) is allowed to shrink with time. For any $\epsilon>0$ and any irrational~$x$ the rotation $qx$ visits the `shrinking target' 
$(\alpha -\epsilon/q, \, \alpha +\epsilon/q) \subset[0,1)$ infinitely often for almost every $\alpha$. On the other hand, Theorem III in Chapter III of Cassels' seminal book~\cite{Cas} shows that (\ref{eqn:Kim}) does not hold for every irrational $x$ and every real $\alpha$:  For every $\epsilon>0$ there is an irrational $x$ such that the pair of inequalities 
$$
\infabs{q} \leq Q, \quad \norm{qx-\alpha} < \frac \epsilon Q
$$
is insoluble for infinitely many values of $Q$. In this sense Kim's result is best possible.

Inspired by this and statement (\ref{eqn:Kim}), investigation into the complementary Lebesgue null set 
\begin{equation*}
	\Bad_{x} \: = \: \left\{\alpha \in [0,1): \, \inf_{q \in \mathbb{N}} q\norm{qx - \alpha } > 0 \right\},
\end{equation*}
quickly followed. 
In 2010 it was shown by Bugeaud et al \cite{BHKV} that this set, and its natural generalisation to higher dimensions, is of maximal Hausdorff dimension. Shortly afterwards, Tseng \cite{Tseng}
demonstrated that $\Bad_{x}$ enjoys the stronger property of being winning 
(in the sense of Schmidt\footnote[1]{We refer the reader to \cite{Schm1} and \cite{Schm2} for all 
necessary definitions and results on winning sets. We only recall here that winning sets in $\mathbb R^n$ necessarily have maximal Hausdorff dimension, and that countable intersections of winning sets are again winning.}) for all real numbers $x$.

In this note we consider the following collection of $mn$-dimensional generalisations of the set $\Bad_{x}$, which allow for the rate of approximation in each coordinate to be assigned  a `weight'. 
Let $x_{ji} $ be real numbers  (for $1 \leq i \leq m, 1 \leq j \leq n$) and let
\begin{equation*}
	L_j(\mathbf{q}) \, = \, \sum_{i=1}^m q_ix_{ji} \quad \quad \quad \quad(1 \leq j \leq n)
\end{equation*}
be the related system of $n$ homogeneous linear forms in the variables $q_1, \ldots, q_m$. 
Denote by $\mathcal{L}$ the $n \times m$ real matrix corresponding to the 
real numbers $x_{ji}$ and by $\mat_{n \times m}(\mathbb{R})$ the set of all such matrices.
Then, for any $n$-tuple of real numbers $\mathbf{k}=\left\{k_1,\ldots, k_n \right\}$ such that
\begin{equation}
	\label{eqn:Krestrictions}
	k_j \: > \: 0 \: \: \: \:(1 \leq j \leq n) \quad \text{ and } \quad \sum_{j=1}^n k_j \, = \, 1,
\end{equation}
 define the set 
\begin{equation*}
	\Bad_{\mathcal{L}}(\mathbf{k}, n, m): \, = \, \left\{ \boldsymbol\alpha \in [0,1)^n: \, \inf_{\mathbf{q} \in \mathbb{Z}_{\neq \mathbf{0}}^m}
	\max_{1 \leq j \leq n} \left(\infabs{\mathbf{q}}^{mk_j}\norm{L_j(\mathbf{q}) - \alpha_j}\right) > 0 \right\}.
\end{equation*}
Here, $\infabs{\, . \, }$ denotes the supremum norm and $\boldsymbol\alpha =(\alpha_1, \ldots \alpha_n     )$. For brevity, we will simply write $\Bad_{\mathcal{L}}(n, m)$ for the standard unweighted case `$k_1=\cdots=k_n=1/n$'. 

Recently, Einsiedler \& Tseng \cite{ET} extended the results of \cite{BHKV} and \cite{Tseng} to show amongst other related results that the set $\Bad_{\mathcal{L}}(n, m)$ is winning for any matrix 
$\mathcal{L} \in \mat_{n \times m}(\mathbb{R})$ (see also \cite{Klein} and \cite{Mosh}). However, it appears their method cannot be extended to the weighted setting of the sets $\Bad_{\mathcal{L}}(\mathbf{k}, n, m)$.

Schmidt was the first to consider weighted variants of the badly approximable numbers. In \cite{Schm3}, he introduced sets of the form 
\begin{equation*}
	\Bad(i, j) \: = \: \left\{(x_1, x_2) \in [0,1)^2: \, \inf_{q \in \mathbb{Z}_{\neq 0}}
	\max\left\{ \infabs{q}^{i}\norm{qx_1}, \,\infabs{q}^{j}\norm{qx_2}\right\}  > 0 \right\},
\end{equation*}
for real numbers $i, j > 0$ satisfying  $i+j=1$. Whilst a metric theorem of Khintchine~\cite{Khi24} implies that these sets are of Lebesgue measure zero, Schmidt noted that 
each set is certainly non-empty. Much later, building on the earlier work of 
Davenport~\cite{Dav}, it was proven by Pollington \& Velani \cite{PV} that each set $\Bad(i, j)$ is always of maximal 
Hausdorff dimension. Remarkably, Badziahin, Pollington \& Velani \cite{BPV} subsequently solved a famous conjecture made by 
Schmidt in \cite{Schm3} stating that the intersection of any two of the distinct sets $	\Bad(i, j)$ is non-empty. Moreover, they proved a general result implying that any finite collection of pairs $(i_t, j_t)$ of strictly positive real numbers satisfying $i_t +j_t=1$ (for $1 \leq t \leq k$) the intersection
$$
\bigcap_{t=1}^k \, \Bad(i_t, j_t)
$$
is of maximal Hausdorff dimension. Under a certain mild technical condition, their statement actually extended to countable intersections. In 2013, An \cite{An2} was surprisingly able to improve both results by demonstrating that each set $	\Bad(i, j)$ is winning.

Inspired by these developments, and those of \cite{BHKV}, the following statement was proven in \cite{Harrap}.
\begin{thmharrap}
	For any real  $i, j > 0$ satisfying  $i+j=1$ and any $\mathbf{x} \in \Bad(i, j)$,
	the Lebesgue null set
	\begin{equation*}
		\Bad_{\mathbf{x}}(i, j) \, = \, \left\{(\alpha_1, \alpha_2) \in [0,1)^2: \, 
		\inf_{q \in \mathbb{Z}_{\neq 0}} \max\left\{ \infabs{q}^{i}\norm{qx_1-\alpha_1}, 	
		\,\infabs{q}^{j}\norm{qx_2-\alpha_2}\right\}  > 0 \right\}
	\end{equation*}
	is of full Hausdorff dimension. 
\end{thmharrap}
The restriction that $\mathbf{x} \in \Bad(i, j)$ was later \cite{Mosh2} safely removed from the statement of Theorem A when $i=2/3$ and $j=1/3$. 

The purpose of this note 
is to extend this result to the full setting of $n$ linear forms in $m$ variables and to establish a statement concerning the intersection of such sets. To do this we require to define one final badly approximable set, 
a natural higher dimensional generalisation of $\Bad(i, j)$. For any $n$-tuple of real numbers 
$\mathbf{k}$ satisfying~(\ref{eqn:Krestrictions}) let
\begin{equation*}
	\Bad(\mathbf{k}, n, m) \, = \, \left\{ \mathcal{L} \in \mat_{n \times m}(\mathbb{R}): \, \inf_{\mathbf{q} 
	\in \mathbb{Z}_{\neq \mathbf{0}}^m}
	\max_{1 \leq j \leq n} \left(\infabs{\mathbf{q}}^{mk_j}\norm{L_j(\mathbf{q})}\right) > 0 \right\}.
\end{equation*}
This set is also known to have zero Lebesgue measure and full Hausdorff dimension~\cite{KW}. It was shown by An \cite{An} that the countable intersection of sets of this form has maximal Hausdorff dimension in the case `$m=1, n=2$'.

\subsection{Statement of Results}
\label{sec:TheMainResults}

We prove the following strengthening of Theorem A.

\begin{thm}
	\label{thm:1}
	For any $n$-tuple $\mathbf{k}$ satisfying 
	(\ref{eqn:Krestrictions}) and any matrix 
	$\mathcal{L}\in \Bad(\mathbf{k}, n, m)$ the set $\Bad_{\mathcal{L}}(\mathbf{k}, n, m)$ is 
	$1/2$ winning.	
\end{thm}

We prove Theorem \ref{thm:1} by adapting the proof of Theorem X (Chapter 5) of Cassels' book \cite{Cas}.  
In short, his theorem implies that the set $\Bad_{\mathcal{L}}(n, m)$ is non-empty. We note that removing the assumption that $\mathcal{L}\in \Bad(\mathbf{k}, n, m)$, whilst desirable, does not seem possible using the methods presented here. 
However, our result does give a (partial) affirmative answer to the analogue of Schmidt's conjecture in the setting of irrational rotations of the circle.
\begin{cor}
Fix any arbitrary sequence $\left\{\mathbf{k}_t    \right\}_{t=1}^\infty$ of $n$-tuples of real numbers $\mathbf{k}_t=\left\{k^{(t)}_1,\ldots, k^{(t)}_n \right\}$ satisfying for every $t \in \mathbb{N}$ the conditions
 \begin{equation*}
k^{(t)}_j \: > \: 0 \: \: \: \:(1 \leq j \leq n) \quad \text{ and } \quad \sum_{j=1}^n k^{(t)}_j \, = \, 1.
\end{equation*} 
Then, for any matrix $\mathcal{L}\in \Bad(\mathbf{k}, n, m)$ the intersection 
$$
\bigcap_{t=1}^\infty \, \Bad_{\mathcal{L}}(\mathbf{k}, n, m)
$$
is of maximal Hausdorff dimension $mn$.
\end{cor}
For completion, we mention the following trivial consequence of 
Theorem~\ref{thm:1} in the more familiar two dimensional setting.
\begin{cor}
	\label{cor:1}
	For any real numbers $i, j > 0$ satisfying  $i+j=1$ and any vector $\mathbf{x}\in 
	\Bad(i, j)$ the set $	\Bad_{\mathbf{x}}(i, j)$ is $1/2$ winning.
\end{cor}

\section{Proof of Theorem \ref{thm:1}}
\label{sec:ProofOfTheoremRefThm1}
For simplicity we assume throughout that the group
$G=\mathcal{L}^{T}\mathbb{Z}^n+\mathbb{Z}^m$ has rank $n+m$. This is because Kronecker's Theorem (see \cite{Kron}) then asserts 
that the dual subgroup $\Gamma= \mathcal{L}\mathbb{Z}^m+\mathbb{Z}^n$ is dense in $\mathbb{R}^n$. In the degenerate case when the rank of $G$ is strictly less than $n+m$ it is easily verified that $\left\{\mathcal{L}\mathbf{q}: \, \mathbf{q} \in \mathbb{Z}^m \right\}$ is restricted to at most a 
countable collection $H$ of parallel, positively separated, hyperplanes in $\mathbb{R}^n$. We therefore have 
$\mathbb{R}^n \setminus H = \Bad_{\mathcal{L}}(\mathbf{k}, n, m)$, from which it is easily deduced that $\Bad_{\mathcal{L}}(\mathbf{k}, n, m)$ 
is winning.

In what follows  	
\begin{equation*}
	M_i(\mathbf{u} ) \, = \, \sum_{j=1}^n u_jx_{ji} \quad \quad \quad \quad(1 \leq i \leq m)
\end{equation*}
denotes the transposed set of $m$ homogeneous linear forms in the variables $u_1, \ldots, u_n$ 
corresponding to the matrix $\mathcal{M}=\mathcal{L}^T$ (the dual forms to $L_j$).
Choose a matrix $\mathcal{L} \in \Bad(\mathbf{k}, n, m)$ and assume without loss of generality that we have $k_1 = \max_{\, 1\leq j\leq n}{ k_j}$. 
We begin by utilising the following lemma, which allows us to switch between the matrices in
$\Bad(\mathbf{k}, n, m)$ and the related `dual' set.
The lemma follows from a general transference theorem which can be found in Chapter V of Cassels' book~\cite{Cas}.
\begin{lem}
	\label{lem:duel}   
	Let $\Bad^{\ast}(\mathbf{k}, m, n)$ be the set of matrices 
	$\mathcal{M} \in \mat_{m \times n}(\mathbb{R})$ such that
	\begin{equation*}
		\inf_{\mathbf{u} \in \mathbb{Z}_{\neq \mathbf{0}}^n}
		\max_{1 \leq i \leq m} \left(\max_{1 \leq j \leq n}
		\left(\infabs{\alpha_i}^{1/(mk_j)}\right)\norm{M_i(\mathbf{u} )} \right) > 0.
	\end{equation*}
	Then,
	\begin{equation*}
		\mathcal{M} \: \in \: \Bad^{\ast}(\mathbf{k}, m, n) \quad \Longleftrightarrow 
\quad \mathcal{M}^T = \mathcal{L} \: \in \: \Bad(\mathbf{k}, n, m).
	\end{equation*}
\end{lem}

For any $T \geq1$ and any $(n+1)$ strictly positive real numbers $\beta_1, \ldots, \beta_{n+1}$ define a set
\begin{eqnarray*}
	\Pi_T(\beta_1, \ldots, \beta_{n+1}) & = & \left\{(\mathbf{u}, \mathbf{v}) \in 
	\mathbb{R}^n \times \mathbb{R}^m:
	\, \infabs{u_j} \leq \beta_j\,T^{\, mk_j} \, \, \,  (1 \leq j \leq n)  \right. \\
	& & \left. \quad \quad \quad \quad \quad \quad \quad \quad 
	 \quad \text{ and } \max_{1\leq i \leq m} \infabs{M_i(\mathbf{u})-v_i} \leq \beta_{n+1}T^{-1} \right\}. \
\end{eqnarray*}
For ease of notation we will hereafter consider sets of this type as subsets of $\mathbb{R}^{n+m}$, 
the origin of which will be denoted $\mathbf{0}$. Now, since $\mathcal{L} \in \Bad(\mathbf{k}, n, m)$, 
Lemma~\ref{lem:duel} immediately implies there exists a constant $\gamma=\gamma(\mathcal{L}) \in (0,1)$ such that
\begin{equation*}
	\Pi_T(1, \ldots, 1, \gamma) \, \cap \, \mathbb{Z}^{n+m}  \: = \: 
	\left\{ \mathbf{0} \right\}.
\end{equation*}
However, the set $\Pi_T(\gamma^{-m}, 1 \ldots, 1, \gamma)$ is a convex, symmetric, closed, 
bounded region in space whose volume is given by 
\begin{equation*}
	2\gamma^{-m}T^{\, mk_1} \cdot \prod_{j=2}^n 2\, T^{\, mk_j} \cdot 2^m\gamma^mT^{-1} \, = \, 2^{n+m}.
\end{equation*}
Therefore, by Minkowski's Convex Body Theorem (see Appendix B of \cite{Cas}) we have that
\begin{equation*}
	\Pi_T(\gamma^{-m}, 1 \ldots, 1, \gamma) \, \cap \, \mathbb{Z}^{n+m}  \: \neq \: 
	\left\{ \mathbf{0} \right\}.
\end{equation*}
This means for any $T \geq 1$ there exists at least one integer vector 
$\mathbf{z}= (\mathbf{u}, \mathbf{v})\in \mathbb{Z}^{n+m}$ such that 
\begin{equation*}
	\mathbf{z} \: \: \in \: \: \Pi_T(\gamma^{-m}, 1 \ldots, 1, \gamma) \: \setminus\: \Pi_T(1, \ldots, 1, \gamma).
\end{equation*}
Choose such an integer vector with the smallest possible first coordinate $u_1\geq1$ for which 
$\max_{\, 1\leq i \leq m} \infabs{M_i(\mathbf{u})-v_i}$ attains its minimal value. Denote this vector by
\begin{equation*}
	\mathbf{z}(T): \: = \: (\mathbf{u}(T), \mathbf{v}(T)): \: = \: (u_1(T), \ldots u_n(T), v_1(T), \ldots, v_m(T)).
\end{equation*}
Also, let
\begin{equation*}
	\phi(T) \: = \: \max_{1\leq i \leq m} \norm{M_i(\mathbf{u}(T))} \: = 
	\: \max_{1\leq i \leq m} \infabs{M_i(\mathbf{u}(T))-v_i(T)}
\end{equation*}
be the minimal value taken.  
Note that the rank assumption imposed on $\mathcal{L}$ ensures that the `best approximation vector' $\mathbf{z}(T)$ always exists and is 
unique up to sign change (for similar constructions, 
see \cite{Lag} or Section 2 of~\cite{BL}).

The following set of inequalities will be useful. Since $\mathbf{z}(T) \in \Pi_T(\gamma^{-m}, 1\ldots, 1, \gamma)$
 we have
\begin{equation}
	\label{eqn:useful1}
	\infabs{u_1(T)} \: \leq \: \gamma^{-m}\, T^{\,mk_1} , \quad \quad \quad \infabs{u_j(T)} \: \leq \: T^{\,mk_j} 
	\quad (2 \leq j \leq n) 
\end{equation}
and also
\begin{equation}
	\label{eqn:useful2}
	\phi(T) \: \leq \: \gamma \,T^{-1}.
\end{equation}
Morover, since $\mathbf{z}(T) \notin \Pi_T(1, \ldots, 1, \gamma)$ we know 
\begin{equation}
	\label{eqn:useful3a}
	\infabs{u_1(T)} \: > \: T^{\,mk_1} \quad \text{ and so } \quad 
	\max_{1 \leq j \leq n} \left( \infabs{u_j(T)}^{1/(mk_j)} \right) \: = \: \infabs{u_1(T)}^{1/(mk_1)}.
\end{equation}
Recalling that $\mathcal{M} \in \Bad^{\ast}(\mathbf{k}, m, n)$, we therefore have
\begin{equation}
	\label{eqn:useful3}
	\phi(T) \: \geq \: \gamma \left( \max_{1 \leq j \leq n} \left( \infabs{u_j(T)}^{1/(mk_j)} \right) \right)^{-1} 
	\: = \: \gamma \infabs{u_1(T)}^{-1/(mk_1)} \: \geq \: \gamma^{1+1/k_1} \, T^{\,-1}.
\end{equation}

Next, we prove a lemma regarding the rate of growth of a suitable sequence of the Euclidean norms 
of the integer vectors $\mathbf{u}(T)$ (c.f.\ \cite[Theorem 1.2]{Lag}).
Put $R:=\left\lceil \gamma^{-1/k_1} \right\rceil +1$ and define $T_r=R^r$ (for $r=0, 1, \ldots$). For 
notational convenience let $\mathbf{z}_r = (\mathbf{u}_r, \mathbf{v}_r) = \mathbf{z}(T_r)$ and 
$\phi_r=\phi(T_r)$.
Inequality (\ref{eqn:useful3}) yields that $\phi_r$ is strictly decreasing as
\begin{eqnarray*}
	\phi_r \: \: \geq \: \: \gamma^{1+1/k_1} \, T_r^{\,-1} \: \: = \: \: \gamma^{1+1/k_1}\, R\, T_{r+1}^{\,-1} 
	& \geq & \gamma^{1+1/k_1}\left( \gamma^{-1/k_1} +1 \right) \, T_{r+1}^{\,-1} \\
	& > & \gamma \, T_{r+1}^{\,-1} \\
	& \geq & \phi_{r+1}.\
\end{eqnarray*}
The final inequality follows from (\ref{eqn:useful2}), which also implies 
\begin{equation}
	\label{eqn:useful4}
	\phi_r \: \: \leq \: \: \gamma \, R \, T_{r+1}^{\, -1}.
\end{equation} 
This will be utilised later, as will the observation that $\phi_r \rightarrow 0$ as $r \rightarrow \infty$.

\begin{lem} 
	\label{lem:lag}
	The sequence of vectors $\left\{\mathbf{u}_r\right\}_{r=0}^{\infty}$ can be partitioned into finitely many subsequences in such a way that the Euclidean norms of the vectors of each subsequence form a lacunary sequence.
\end{lem}

\vspace{1cm}

\begin{proof}
		Consider the Euclidean norm $\infabs{ \, . \,}_e$ of each integer vector $\mathbf{u}_r$.
	From (\ref{eqn:useful3a}) we have
	\begin{eqnarray}
		T_r^{\, 2mk_1} \: \: < \: \: \infabs{u_1(T_r)}^2 & \leq & \infabs{\mathbf{u}_r}_e^2 
		\label{eqn:useful5a}\\ 
		& \stackrel{\text{(\ref{eqn:useful1})}}{\leq} & \gamma^{-2m}\, T_r^{\, 2mk_1} + \sum_{j=2}^{n} T_r^{\, 2mk_j}
		\nonumber \\ 
		& < & \gamma^{-2m} \sum_{j=1}^n T_r^{\, 2mk_j} \nonumber \\
		 & \leq &  \gamma^{-2m} n \, T_r^{\, 2mk_1}, \label{eqn:useful5b}\
	\end{eqnarray}
	since we are assuming $k_1 = \max_{\, 1\leq j\leq n}{ k_j}$. Now, choose any natural number $t$ such that 
	$R^{tmk_1} \geq 2n^{1/2} \gamma^{-m}$. Then,
	\begin{equation*}
		\infabs{\mathbf{u}_{r+t}}_e \: \:\stackrel{\text{(\ref{eqn:useful5a})}}{>} \: \: T_{r+t}^{\, mk_1} \: \:
		= \: \: R^{tmk_1} \, T_r^{\, mk_1} \: \: \geq \: \: 2n^{1/2} \gamma^{-m}\, T_r^{\, mk_1} \: \: 
		\stackrel{\text{(\ref{eqn:useful5b})}}{>} \:\: 2\infabs{\mathbf{u}_{r}}_e.
	\end{equation*}
	So, the sequence $\left\{\mathbf{u}_r\right\}_{r=0}^{\infty}$ can be partitioned into a finite collection 
	of subsequences $\left\{\mathbf{u}_{t_0+tr}\right\}_{r=0}^{\infty}$ such that each subsequence is $2$-lacunary;
	that is
	\begin{equation*}
				\infabs{\mathbf{u}_{t_0+t(r+1)}}_e \quad \geq \quad 2 \infabs{\mathbf{u}_{t_0+tr}}_e \quad \forall \, r.
	\end{equation*}
	
\end{proof}
Note that in order to construct the above sequences the assumption that $\mathcal{M} \in \Bad^{\ast}(\mathbf{k}, m, n)$ was imperative. The lemma allows us to utilise the following powerful result, which is taken from~\cite{Mosh}.
\begin{lem}
	\label{lem:mosh}
	If a sequence $\left\{ \mathbf{w}_r\right\}_{r=0}^{\infty}$ of non-zero integral vectors
	 is such that the 
	corresponding sequence of	Euclidean norms is lacunary then the set
	\begin{equation*}
		\left\{ \boldsymbol\alpha \in [0,1)^n: \, \inf_{r} 
		\norm{\mathbf{w}_r \cdot \boldsymbol\alpha } > 0 \right\}
	\end{equation*}
	is $1/2$ winning.
\end{lem}
\begin{cor}
	\label{cor:mosh}
	The set
	\begin{equation*}
		\Bad_{\left\{ \mathbf{u}_r \right\}} \:\: = \: \: \left\{ \boldsymbol\alpha \in [0,1)^n: \, \inf_{r} 
		\norm{ \mathbf{u}_r  \cdot \boldsymbol\alpha } > 0 \right\}
	\end{equation*}
	is $1/2$ winning.
\end{cor}
Corollary \ref{cor:mosh} follows from Lemma \ref{cor:mosh} by the observation that
\begin{equation*}
\Bad_{\left\{ \mathbf{u}_r \right\}} \: \: = \: \: \bigcap_{t_0=0}^{t-1} \left\{ \boldsymbol\alpha \in \mathbb{R}^n: \, \inf_{r} 
\norm{\mathbf{u}_{t_0+tr}  \cdot \boldsymbol\alpha } > 0 \right\}.
\end{equation*}

The set $\Bad_{\left\{ \mathbf{u}_r \right\}}$ was first shown in~\cite{BHKV} to have full Hausdorff dimension for any sequence $\left\{ \mathbf{u}_r\right\}_{r=0}^{\infty}$ of non-zero integral vectors whose Euclidean norms form a lacunary sequence. We remark that in the case where the matrix   $\mathcal{M}$ is not assumed to be chosen from $\Bad^{\ast}(\mathbf{k}, m, n)$  one may partition the sequence $\left\{\mathbf{u}_r\right\}_{r=0}^{\infty}$ into a finite collection of subsequences $\left\{\mathbf{u}_{t_0+tr}\right\}_{r=0}^{\infty}$ such that each subsequence is lacunary with respect to the norm in (\ref{eqn:useful3a}); that is, the norm given by
$$\infabs{\mathbf{x}}_{\mathbf{k}}:=\max_{1 \leq j \leq n} \left( \infabs{x_j}^{1/(mk_j)} \right) .$$ 
However, showing that the analogue of Lemma \ref{cor:mosh} holds in this case does not seem straightforward. Moreover, it is the belief of the authors that the corresponding set may not in fact be winning.

We are now ready to prove Theorem \ref{thm:1}. Choose $\boldsymbol\alpha \in \Bad_{\left\{ \mathbf{u}_r \right\}}$ and assume
\begin{equation*}
	\inf_{r} \norm{\mathbf{u}_r  \cdot \boldsymbol\alpha}  \: \geq \: \epsilon  \: >  \: 0.
\end{equation*}
For any $\mathbf{q} \in \mathbb{Z}_{\neq \mathbf{0}}^m$, the trivial equality
\begin{equation*}
	\mathbf{u}_r  \cdot \boldsymbol\alpha\, = \,  \sum_{i=1}^{m}q_iM_i(\mathbf{u}_r) \, - \, 
	\sum_{j=1}^n(L_j(\mathbf{q})-\alpha_j)u_{j}(T_r), 
\end{equation*}
in conjunction with the triangle inequality yields that
\begin{eqnarray}
	0 < \epsilon <	\norm{\mathbf{u}_r  \cdot \boldsymbol\alpha} & \leq &  
	m\max_{1\leq i \leq m} \left(  \norm{M_i(\mathbf{u}_r)}\infabs{q_i} \right) \: \: + \: \: n \max_{1\leq j \leq n}
	\left(\norm{(L_j(\mathbf{q})-\alpha_j)}\infabs{u_j(T_r)} \right)\nonumber \\
	& \leq & m\phi_r\infabs{\mathbf{q}}  \: \: + \: \:
	n\max_{1\leq j \leq n}\left(\norm{(L_j(\mathbf{q})-\alpha_j)}\infabs{u_j(T_r)} \right).     \label{eqn:yuvnorm} \
\end{eqnarray}
Here, we have employed the fact that $\norm{a\mathbf{z}} \leq \infabs{a}\norm{\mathbf{z}}$ for all $a\in \mathbb{R}$ 
and all $\mathbf{z}\in \mathbb{R}^k$. 

Since $\phi_r$ is strictly decreasing and $\phi_r \rightarrow 0$ as $r \rightarrow \infty$ we are free to choose 
$r$ in such a way that 
\begin{equation}
	\label{eqn:defr}
	\phi_r \, < \, \frac{\epsilon}{2m \, \infabs{\mathbf{q}}} \, \leq \phi_{r-1},
\end{equation}
whereby inequality (\ref{eqn:yuvnorm}) yields
\begin{equation*}
	\max_{1\leq j \leq n}\left(\norm{(L_j(\mathbf{q})-\alpha_j)}\infabs{u_j(T_r)} \right)\: \geq \: \epsilon/2n.
\end{equation*}
Finally, notice that combining (\ref{eqn:useful4}) with (\ref{eqn:defr}) implies 
\[T_r \, \leq \, 2m\epsilon^{-1} 
\gamma R  \infabs{\mathbf{q}},\]
and so we have
\begin{equation*}
	\infabs{u_1(T_r)} \, \stackrel{\text{(\ref{eqn:useful1})}}{\leq} \, \gamma^{-m} \, T_r^{\, mk_1} \, 
	\leq \, (2mR)^{mk_1}\gamma^{m(k_1-1)}\epsilon^{-mk_1}  \infabs{\mathbf{q}},
\end{equation*}
and similarly (for $2\leq j \leq n$)
\begin{equation*}
	\infabs{u_j(T_r)} \, \leq \, (2mR\gamma)^{mk_1}\epsilon^{-mk_1} \infabs{\mathbf{q}}.
\end{equation*}
Therefore, 
\begin{equation*}
	\max_{1 \leq j \leq n} \left(\norm{L_j(\mathbf{q}) - \alpha_j}\infabs{\mathbf{q}}^{mk_j}\right) 
	\: \: \geq \: \: \kappa,
\end{equation*}
for some constant $\kappa>0$.
Since the choice of vector $\mathbf{q}$ was arbitrary we have shown that $\boldsymbol\alpha \in \Bad_{\mathcal{L}}(\mathbf{k}, n, m)$, and in particular that $\Bad_{\left\{ \mathbf{u}_r \right\}} \subseteq \Bad_{\mathcal{L}}(\mathbf{k}, n, m)$. In view of Corollary \ref{cor:mosh}, 
the desired conclusion easily follows.

\providecommand{\bysame}{\leavevmode\hbox to3em{\hrulefill}}
\providecommand{\MR}{\relax\ifhmode\unskip\space\fi MR }
\providecommand{\MRhref}[2]{%
  \href{http://www.ams.org/mathscinet-getitem?mr=#1}{#2}
} \providecommand{\href}[2]{#2}

\end{document}